\documentclass[amssymb,amstex,amsmath,10pt]{article}
\usepackage[margin=1.4in]{geometry}
\usepackage{amsmath,amssymb,latexsym,amsfonts,url}
\usepackage{graphicx}
\usepackage{amsthm}
\usepackage{enumitem}
\usepackage{amssymb}
\usepackage{color}
\usepackage{epstopdf}
\renewcommand{\epsilon}{\varepsilon}
\theoremstyle{plain}
\newtheorem{thm}{Theorem}[section]
\newtheorem{prop}[thm]{Proposition}

\newtheorem{cor}[thm]{Corollary}
\newtheorem{lem}[thm]{Lemma}

\newtheorem{exam}[thm]{Example}

\newtheorem*{theorem*}{Theorem}
\newtheorem*{proposition*}{Proposition}
\theoremstyle{definition}

\theoremstyle{remark}

\numberwithin{equation}{section}

\newcommand{\tr}{{\rm{tr}}}

\def\({\left(}
\def\){\right)}

\newcommand{\be}{\begin{eqnarray}}
\newcommand{\ee}{\end{eqnarray}}
\newcommand{\bea}{\begin{eqnarray*}}
\newcommand{\eea}{\end{eqnarray*}}

\newcommand{\n}{\nabla}
\def\e{\epsilon}

\def\o{\omega}
\def\t{\theta}
\def\tr{\Delta}
\def\a{\alpha}
\def\vp{\varphi}

\def\d{\delta}

\def\dis{\displaystyle}

\def\vp{\varphi}

\begin{document}
\begin{title}
{On the geometry of Einstein-type manifolds with some structural conditions}
\end{title}
\begin{author}
{Gabjin Yun \and Seungsu Hwang$^*$}
\end{author}

\date{\today}

\maketitle

\begin{abstract}
\noindent
 In this paper, we investigate the geometry of   Einstein-type equation on  a Riemannian manifold,  unifying various particular geometric structures recently studied in the literature, such as critical point equation  and vacuum static equation. We show various rigidity results of Einstein-type manifolds
  under assumptions of several  curvature conditions.

\vspace{.12in}


\end{abstract}

\setlength{\baselineskip}{15pt}

\section{Introduction}

In the  last two decades since a solution of  geometrization conjecture due to Perelman \cite{per01, per02, per03}, there has been  increasing interest in the study of Riemannian manifolds endowed with metrics satisfying some geometric structural equations, possibly involving curvatures and some globally defined smooth functions called potential functions. Typical examples are  gradient Ricci solitons arisen as self-similar solutions of the Ricci flow (cf. \cite{cao}, \cite{e-l-m}), a
critical point equation  arisen as  an Euler-Lagrange equation of the total scalar curvature functional restricted to the metrics of constant scalar curvature with unit volume on a compact smooth manifold (cf. \cite{Be}, p.128), and vacuum static equations handled in mathematical physics or general relativity (cf. \cite{h-e}).

In this paper we consider an $n$-dimensional smooth Riemannian manifold  $(M,g)$ with
 $n\geq 3$ which admits  smooth functions $f$ and $h$ to the system of equations
\be 
f {\rm Ric}= Ddf +h g,\label{eq1}
\ee
which called the Einstein-type equation. Here, ${\rm Ric}$ is the Ricci curvature of $(M, g)$
 and $Ddf$ denotes  the Hessian of $f$.
We say that the quadruple $(M,g, f, h)$ is called an  Einstein-type manifold if $(M, g)$ is a smooth Riemannian manifold and the pair $(f,h)$ satisfies the Einstein-type equation (\ref{eq1}).

The notion of Einstein-type equations is widely handled in several papers. For example, 
Qing and Yuan \cite{QY}  considered Riemannian manifolds satisfying Einstein-type equation, and 
obtained rigidity results on critical point equation and vacuum static spaces under completely divergence-freeness of the Cotton tensor. Catino et al. \cite{c-m-m-r} provided a more general Einstein-type equation and classified it under the Bach-flat condition. Recently, Leandro \cite{lea} also studied the Einstein-type equation, and showed harmonicity of Weyl curvature tensor under completely divergence-freeness of Weyl curvature and zero radial Weyl curvature condition. 

\vskip .75pc

\noindent
{\bf Notations: }  Hereafter, for convenience and simplicity, we denote the Ricci curvature 
${\rm Ric}$ just by $r$ if there is no ambiguity.
We also denote by $s$ the scalar curvature of $(M, g)$ and, if necessary, we
use the notation $\langle \,\, ,\,\, \rangle$ for metric $g$ or inner product induced by $g$ on tensor spaces.

\vskip .75pc

First of all, by taking the trace of (\ref{eq1}), we trivially obtain
\be 
\tr f=sf-nh. \label{eq7227}
\ee
If $f \equiv 0$, then $h$ must be zero and  the Einstein-type equation becomes trivial having no informations. Thus,  we assume $f$ is not identically zero whenever we consider the Einstein-type equation. It is also easy to see that  $g$ is Einstein if $f$ is constant, and so, in particular, space forms are Einstein-type manifolds. We sometimes call an Einstein-type manifold $(M, g,f,h)$ is {\it trivial}
 if $f$ is a non-zero constant.

As mentioned above, some well known structural equations are directly related to Einstein-type manifolds. For example, if $h=0$, it reduces to the static vacuum Einstein equation (cf. \cite{and, h-c-y}).
 If $h=\frac s{n-1}f$, then we have $s_g'^*(f)=0$, or the vacuum static equation (cf. \cite{amb, 
 f-m, kob, laf}), where $s_g'^*$ is the $L^2$-adjoint operator of the linearization $s_g'$ of the scalar curvature given by
$$
s_g'^*(f)= -(\tr f)g +Ddf-fr.
$$
If $h=\frac s{n-1}f+\frac {\kappa}{n-1}$ for some constant $\kappa$, we have $s_g'^*(f)=\kappa g$, or the $V$-static equation. Also, if $h$ is constant with $s=0$, then $(M,g,f)$ satisfies the $V$-static equation \cite{c-e-m}, since
$$ s_g'^*(f)= (n-1)h g.$$
If $h=\frac s{n-1}f-\frac s{n(n-1)}$ with $f=1+\varphi$, we have $s_g'^*(f)= {\mathring {\rm Ric}}$, called the critical point equation, where ${\mathring {\rm Ric}}$ denotes the trace-less Ricci tensor. Finally, if $h=\frac {s-\rho-\mu}{n-1}f$, we have the static perfect fluid equation \cite{k-o}. By considering $-f$ instead of $f$, we always assume that $h \ge 0$ when $h$ is constant.

In Section 2, we see that if the scalar curvature is vanishing, then the function $h$ must be constant, and conversely, if $h$ is constant and $f$ is not constant, then the scalar curvature $s$ must  vanish. 
In fact, we see that if both $h$ and the scalar curvature $s$ are constants, then either
$(M, g)$ is Einstein or $s =0$.
Moreover, when the function $h$ is vanishing, the Einstein-type equation becomes, in fact, a static vacuum Einstein equation (cf. \cite{wal}) with zero scalar curvature if $M$ is compact.
The first result  is a rigidity theorem when $h=0$.

\begin{thm} \label{thm1} 
Let $(M,g,f)$ be an $n$-dimensional compact Einstein-type manifold satisfying
\be
fr = Ddf.\label{eqn2022-3-5-1}
\ee
Then,  $(M,g)$ is Ricci flat with constant $f$. 
\end{thm}
We would like to mention  a remark on noncompact manifolds satisfying (\ref{eqn2022-3-5-1}).
In addition to (\ref{eqn2022-3-5-1}), if we assume $\Delta f = 0$,  Theorem~\ref{thm1}  is still true for a complete noncompact Riemannian manifold $(M, g)$ (cf. \cite{lic, and0, and}). 
Of course, if we do not assume $\Delta f = 0$, Theorem~\ref{thm1} does not hold anymore (see Section 2). In fact, in case of non-compact Riemannian manifolds which are not complete, 
there are much more complicated geometric structures in static vacuum Einstein manifolds (cf. \cite{and}).

In case when a compact Einstein-type manifold has positive scalar curvature, by applying the maximum principle, we can show  the following a gap-type theorem.

\begin{thm}\label{thm4}
Let $(M,g,f,h)$ be a compact Einstein-type manifold with positive scalar curvature with constant $h$.
 If $ {\min_M  s} \geq nh$, then $(M,g)$ is Einstein.
\end{thm}

When the Ricci curvature is nonpositive, it is easy to see from the Bochner-Weitzenb\"ock formula that the maximum principle holds for the function $|\n f|^2$. 
In particular, we  have the following result.

\begin{thm} \label{thm930} 
Let $(M^n,g,f, h)$ be a complete noncompact Einstein-type manifold with nonpositive Ricci curvature satisfying (\ref{eq1}) with constant $h$.
If $f$ satisfies 
$$
\int_{M}|\nabla f|^2 < \infty,
$$
then $f$ is constant and $(M, g)$ is Einstein.
\end{thm}

As an immediate consequence, we have the following result.

\begin{cor}\label{cor930}
Let $(M^n,g,f,h)$ be a compact Einstein-type manifold with nonpositive Ricci curvature satisfying (\ref{eq1}) with constant $h$.
Then $(M, g)$ is Einstein with constant $f$.
\end{cor}

Next, we consider rigidity question on Einstein-type manifolds with locally conformally flat structure.
Since an Einstein  manifold having locally conformally flat structure has constant sectional curvature,
it is natural to consider Einstein-type manifolds with locally conformally flat metrics.
 When  an Einstein-type manifold $(M,g, f, h)$ is locally conformally flat,  we have the following result.
 
\begin{thm}\label{thm129}
Let $(M,g,f,h)$ be a locally conformally flat Einstein-type manifold satisfying (\ref{eq1}) with constant $h$ and constant scalar curvature.  If $f$ is a proper map, then $(M, g)$ is Einstein.
\end{thm}

For an Einstein-type manifold satisfying conditions in Theorem~\ref{thm129},
we can show that, around any regular point of $f$, the manifold is locally a warped product of an interval and a $(n-1)$-dimensional space form. Related to conformally flat Einstein-type manifolds, we would like to mention a result due to Kobayashi and Obata in \cite{k-o}. They proved that if a warped product manifold $\widetilde M = {\Bbb R} \times M$ is locally conformally flat and the base manifold $(M, g)$   satisfies 
 the Einstein-type equation by a positive potential function, then $(M, g)$ is of  constant curvature.
 
We say an Einstein-type Riemannian manifold $(M, g)$ satisfying (\ref{eq1})
has  radially flat Weyl curvature if $\tilde{i}_{\nabla f}{\mathcal W}=0$. Here
$\tilde{i}_{\nabla f}$ denotes the interior product to the last component defined by
$$
\tilde{i}_{\nabla f}{\mathcal W}(X, Y, Z) = {\mathcal W}(X, Y, Z, \n f).
$$
Note that if the dimension of $M$ is four, an Einstein-type manifold having radially flat Weyl curvature is locally conformally flat (see the proof of Lemma 4.3 in \cite{cc}). Thus, as an immediate consequence of the above theorem, we have the following result. 

\begin{cor}
Let $(M,g,f,h)$ be a $4$-dimensional Einstein-type manifold with radially flat Weyl curvature
satisfying (\ref{eq1}) with constant $h$ and  constant scalar curvature.  If $f$ is a proper map, then $(M, g)$ is Einstein.
\end{cor}

Since $h$ must be constant if the scalar curvature is vanishing for an Einstein-type manifold
$(M, g,f,h)$ satisfying (\ref{eq1}),  we have the following result.

\begin{cor} 
Let $(M,g,f,h)$ be a locally conformally flat Einstein-type manifold with vanishing scalar curvature. If $f$ is a proper map, then $(M, g)$ is flat.
\end{cor}

The paper is organized as follows. In Section 2, we derive basic facts on Einstein-type manifolds in general and Einstein-type manifolds with positive scalar curvature, and prove Theorem~\ref{thm4}.
In section 3, we handle Einstein-type manifolds with constant $h$ or $h=0$ and having nonpositive Ricci curvature, and prove Theorem~\ref{thm1} and  Theorem~\ref{thm930}.  In Section 4, we study Einstein-type manifolds with zero Cotton tensor (see Section 4 for the definition) or locally conformally flat structure.


\section{Basic Properties}

In this section, we shall find basic properties of the scalar curvature of Einstein-type manifolds
 satisfying (\ref{eq1}). First of all, since $\d Ddf = -r(\n f, \cdot) - d\Delta f$ and $\Delta f = fs - nh$, by taking the divergence $\d$ of (\ref{eq1}), we obtain
\be
dh= \frac 1{2(n-1)}\left( fds +2sdf\right). \label{eq83}
\ee
This identity shows that if $s = 0$, then $h$ must be constant. The converse is also true if both $s$ and $h$ are constants and $f$ is not constant.
From (\ref{eq83}), it is easy to see that the following equalities hold in general.

\begin{lem} \label{cor722} On an Einstein-type manifold $(M,g,f,h)$ we have
$$Ddh=\frac 1{2(n-1)}df\otimes ds +\frac 1{2(n-1)}f Dds +\frac 1{n-1}ds\otimes df +\frac s{n-1}(fr -hg).
$$
In particular,
$$\tr h = \frac 3{2(n-1)}\langle \n s, \n f\rangle +\frac 1{2(n-1)}f\tr s +\frac {s^2f}{n-1}  -\frac {nhs}{n-1}.
$$
\end{lem}


A direct observation from (\ref{eq83}) is the following.

\begin{lem}\label{lem2022-1-26-1}
Let $(M^n, g,f, h)$ be an Einstein-type manifold satisfying (\ref{eq1}).
If $h$ is constant, then the function $f^2 s$ must be constant.
\end{lem}
\begin{proof}
From (\ref{eq83}), we have
$$
d(f^2 s) = f(fds + 2sdf) = 0.
$$
\end{proof}

\begin{prop}\label{prop2021-11-11-32}
Let $(M, g, f)$ be a compact Einstein-type manifold   satisfying
(\ref{eq1}) for  a constant $h $. Then
$$
 \int_M f(nh  - fs) dv_g  \ge 0.
$$
The equality holds if and only if $(M, g)$ is Einstein.
\end{prop}
\begin{proof}
 Multiplying both sides of $\Delta f = fs - nh$ by $f$ and integrating it over $M$, we obtain
$$
\int_M f(nh  - fs) dv_g = - \int_M f \Delta f dv_g = \int_M |\n f|^2 dv_g \ge0.
$$
The equality holds if and only if $f$ is constant and so $(M, g)$ is Einstein. 
\end{proof}

In the rest of  this section,  we  discuss the rigidity of Einstein-type manifolds when $h=0$ and  we prove Theorem~\ref{thm1}.
To do this,  we need the following lemma which shows that a compact Riemannian manifold satisfying the  static vacuum Einstein equation must have nonnegative scalar curvature.

\begin{lem}\label{lem83} 
Let $(M,g,f)$ be an $n$-dimensional compact Einstein-type manifold 
satisfying $fr = Ddf$. Then the scalar curvature is nonnegative.
\end{lem}
\begin{proof} First, we claim that there are no critical points of $f$ on $f^{-1}(0)$ unless $f^{-1}(0)$ is empty. Suppose that  $p\in f^{-1}(0)$ is a critical point of $f$.
Let $\gamma$ be a normal geodesic starting from  $p$ moving toward the inside of 
$M^+=\{x\in M\, \vert\, f(x) >0\}$. Then, by (\ref{eq1}), the function 
 $\varphi(t)=f\circ \gamma (t)$ satisfies
 $$ 
 \varphi''(t)= Dd\varphi(\gamma'(t), \gamma'(t))=\varphi(t)r(\gamma'(t), \gamma'(t)).
 $$
Since $\varphi(0)=0$ and $\varphi'(0)=df_p(\gamma'(0))=0$, by the uniqueness of ODE solution,
 $\varphi$ vanishes identically, implying that $\gamma(t)$ stays in $f^{-1}(0)$, which is a contradiction. 
 
 Since $ s\, df=0$ on the set $f^{-1}(0)$ by (\ref{eq83}), and there are no critical points on 
$f^{-1}(0)$, the scalar curvature vanishes on $f^{-1}(0)$.
Now, by (\ref{eq83}) again, we have
$$ 
\frac 32 \langle \nabla s, \nabla f\rangle +s\tr f +\frac f2 \tr s=0
$$
with $\tr f =sf$, implying that 
\be 
\tr s -\frac {6|\nabla f|^2}{f^2}s=-2s^2\leq 0 \label{eqn2022-1-26-2}
\ee
on the set $f \ne 0$. For a sufficiently small positive real number $\e>0$, let
$$
M^\e =\{x\in X\,:\, f(x) > \e\}.
$$
Applying  the maximum principle (cf. \cite{g-t})  to (\ref{eqn2022-1-26-2}) on the set
$M^\e$, we have
$$ 
\inf_{M^\e} \, s \geq \inf_{\partial M^\e}\, s^{-}, 
$$
where $s^{-} = \min\{s, 0\}$.
By letting $\e \to 0$, we have
$$ 
\inf_{M^0} \, s \geq \inf_{\partial M^0}\, s^{-} =0, 
$$
where  $M^0=\{x\in M\, \vert\, f(x) >0\}$ with $\partial M^0=f^{-1}(0)$.
The equality in right-hand side follows from the above claim. In a similar way, we may argue that 
$$ 
\inf_{M_0} \, s \geq \inf_{\partial M_0}\, s^{-} =0,
$$
where $M_0=\{x\in M\, \vert\, f(x) <0\}$ with $\partial M_0=f^{-1}(0)$.
As a result, we may conclude that $s\geq 0$ on $M$. 

Finally, assume that  $f^{-1}(0)$ is empty. Letting $\min_M s = s(x_0)$, we have
$\Delta s(x_0) = 0$ by Lemma~\ref{lem2022-1-26-1} and (\ref{eqn2022-1-26-2}) and so  $s(x_0) = 0$   from (\ref{eqn2022-1-26-2}) again.
\end{proof}

Now, we are ready to prove Theorem~\ref{thm1}. 

\vspace{.12in}
\noindent
{\bf Proof of Theorem~\ref{thm1}.}
 By Lemma~\ref{lem83}, we have  $s\geq 0$ on $M$, and so, from $\tr f =sf$, we have
$$\tr f\geq 0$$
on the set $M^0$. Since $f \not\equiv 0$, from the maximum principle, we can see
that there is only one case:  $f$ is a nonzero constant and $s = 0$ on $M$. Thus
 $(M, g)$ is Ricci-flat. \hfill $\Box$

\begin{exam}
{\rm
Let $M^3 = (a, b) \times {\Bbb R}^2$ be a smooth $3$-manifold with the metric
\be
g = \frac{dt^2}{f^2} + t^2(d\t^2+\sin^2 \t d\vp^2 ),\quad f: (a,b) \to {\Bbb R}^+ \,\,\,{\rm smooth.}\label{eqn2021-1-7-1}
\ee
Recall that $d\t^2+\sin^2 \t d\vp^2$ is the standard spherical metric on ${\Bbb S}^2$ so that
$dt^2 +t^2(d\t^2+\sin^2 \t d\vp^2)$ is just the flat metric on ${\Bbb R}^3$.
A standard frame is 
$$e_1 = f\frac{\partial}{\partial t}, \quad e_2 = \frac{1}{t}\frac{\partial}{\partial \t},  \quad e_3 = \frac{1}{t\sin\t}\frac{\partial}{\partial \vp}$$ 
and the corresponding coframe is 
$$
\o^1 = \frac{1}{f}dt,  \quad \o^2 = t d\t,\quad \o^3 = t\sin \t d\vp.
$$
Using connection form and curvature form from these, it is easy to compute the Ricci curvature as follows: denoting $R_{ij} = {\rm Ric}(e_i, e_j)$,
\be
R_{11} = - \frac{2f'f}{t}, \quad R_{22} = - \frac{tf'f +f^2 -1}{t^2},\quad 
R_{33} =  - \frac{tf'f +f^2 -1}{t^2} = R_{22}\label{eqn2022-3-5-2}
\ee
and
$$
R_{12} = 0 = R_{13} = R_{23}.
$$

\vspace{.15in}
\noindent
Here we consider two cases:
\begin{itemize}
\item[(1)] $\dis{f(t) = \sqrt{1-\frac{2m}{t}}}$ \,\,\, and \,\,\, $(a,b) = (2m, \infty)$
\item[(2)] $\dis{f(r) = \sqrt{1-\frac{2mt^2}{R^3}}, \,\, R>0}$ \,\,\, and \,\,\, $\dis{(a,b) = \left(-\sqrt{\frac{R^3}{2m}}, \,\, \sqrt{\frac{R^3}{2m}}\, \, \right)}$
\end{itemize}

\vspace{.152in}
\noindent
{\bf Case (1)}: $\dis{f(r) = \sqrt{1-\frac{2m}{t}}}$ \,\,\, and \,\,\, $(a,b) = (2m, \infty)$

\vspace{.15in}
In this case, we have  
$$
f' = \frac{1}{2f}\cdot \left(- \frac{2m}{t}\right)' = \frac{m}{t^2f},\quad f'f = \frac{m}{t^2}.
$$
So, by (\ref{eqn2022-3-5-2}), we obtain
$$
R_{11} = - \frac{2m}{t^3}, \quad R_{22} = R_{33} = \frac{m}{t^3}
$$
and
$$
s = R_{11} + R_{22} + R_{33} = 0.
$$

\vspace{.15in}
\noindent
{\bf Case (2)}: \,\, $\dis{f(t) = \sqrt{1-\frac{2mt^2}{R^3}}}$\,\,\, and \,\,\, $\dis{(a,b) = \left(-\sqrt{\frac{R^3}{2m}}, \,\, \sqrt{\frac{R^3}{2m}}\, \, \right)}$

\vspace{.15in}
In this case, we have $\dis{\frac{f^2 - 1}{t^2} = - \frac{2m}{R^3}}$ and
\be
f' =- \frac{2mt}{R^3 f}, \quad \frac{f'f}{t} = - \frac{2m}{R^3}. \label{eqn2021-1-7-6}
\ee
So, by (\ref{eqn2022-3-5-2})
$$
R_{11} = \frac{4m}{R^3} = R_{22} = R_{33}
$$
and 
$$ s = \frac{12m}{R^3}.
$$
In particular, this is an Einstein manifold.

\vspace{.15in}
\noindent
Now, we compute the Laplacian and Hessian of $f$.
A computation shows that
\bea
\Delta =
f^2 \frac{\partial^2}{\partial t^2} + \left(\frac{2f^2}{t} + f'f\right) \frac{\partial}{\partial t} + \frac{1}{t^2} \frac{\partial^2}{\partial \t^2}
+ \frac{\cos \t}{t^2 \sin \t}\frac{\partial}{\partial \t} + \frac{1}{t^2\sin^2\t}\frac{\partial^2}{\partial \vp^2}.
\eea
For the first case $\dis{f(t) = \sqrt{1-\frac{2m}{t}}}$, we have $\Delta f = 0$ and
the Hessian of $f$ is given by
$$
Ddf(e_1, e_1) = - \frac{2m}{t^3}f, \quad Ddf(e_2, e_2)  =  \frac{m}{t^3}f, \quad
Ddf(e_3, e_3) =  \frac{m}{t^3}f.
$$
Thus, we have
$$
f{\rm Ric} = Ddf.
$$
For the second case   $\dis{f(t) = \sqrt{1-\frac{2mt^2}{R^3}}}$, we have
$$
\Delta f = - \frac{6m}{R^3} f  = - \frac{s}{2} f
$$
and
$$
\frac13 \left(sf - \Delta f\right) = \frac{4m}{R^3}f + \frac{2m}{R^3}f  = \frac{6m}{R^3} f.
$$
The Hessian of $f$ is given as follows.
$$
Ddf(e_1, e_1) = - \frac{2m}{R^3}f, \quad Ddf(e_2, e_2)  = - \frac{2m}{R^3}f, \quad
Ddf(e_3, e_3) = - \frac{2m}{R^3}f.
$$
Hence
\bea
 f {\rm Ric} = Ddf + \frac13 (sf - \Delta f)g.
\eea
\hfill $\Box$}

\end{exam}

\section{Einstein-type manifolds with positive scalar curvature}

In this section, we consider Einstein-type manifolds with positive constant scalar curvature. 
First, the following property shows that if $(M, g, f, h)$ is a compact Einstein-type manifold satisfying (\ref{eq1}) with nonzero constant $h$, and $(M, g)$ has  positive scalar curvature, then we have $h f >0$ on $M$.

\begin{prop}\label{eqn2021-11-1-1}
Let $(M, g, f, h)$ be a compact Einstein-type manifold  with positive scalar curvature satisfying
(\ref{eq1}) for a  nonzero constant $h$. Then
we have $h  f >0$ on the whole $M$, and so
$$
h  \int_M f dv_g >0.
$$
\end{prop}
\begin{proof}
If $f$ is a nonzero constant, then we have $nhf = f^2 s >0$. Now, we assume that $f$ is not constant.
Let 
$$
\min_M f = f(x_0)\quad \mbox{and}\quad \max_M f = f(x_1).
$$
If  $ f(x_1) >0$, then $\Delta  f = fs - nh  \le 0$ at the point $x_1$, and so we have $h>0$.
From $\Delta f(x_0) = fs - nh  \ge 0$ at the point $x_0$, we have
$$
f>0
$$
on the whole $M$. Similarly, if $ f(x_0) <0$, then
$\Delta  f = fs - nh  \ge 0$ at the point $x_0$, and so we have $h <0$. From
$\Delta f(x_1) = fs - nh  \le 0$ at the point $x_1$, we have
$$
f<0
$$
on the whole $M$. Finally, it is easy to see that $h<0$ if $f \le 0$, and $h>0$ if $f \ge0$.
\end{proof}

A  similar proof as in Proposition~\ref{prop2021-11-11-32} shows the following.

\begin{prop}\label{prop2021-11-11-33}
Let $(M, g, f)$ be a compact Einstein-type manifold  with positive scalar curvature satisfying
(\ref{eq1}) for a   positive constant $h$. Then
$$
\int_M s(nh  - fs) dv_g  \le 0.
$$
The equality holds if and only if $(M, g)$ is Einstein.
\end{prop}
\begin{proof}
Multiplying both side of $\Delta f = fs - nh $ by s and integrating it over $M$, we have
$$
\int_M s(nh - fs) dv_g = -\int_M s \Delta f =  \int_M \langle \n s, \n f\rangle
= -2\int_M \frac{s}{f}|\n f|^2 dv_g.
$$
In the last equality, we used the identity $fds = - 2s df$.  
Since  $h >0$  and $h f>0$ by Proposition~\ref{eqn2021-11-1-1}, the proof is complete.
Finally, the equality holds if and only if $f$ is constant and so $(M, g)$ is Einstein. 
\end{proof}

For a nonnegative integer $m$, let us define $\varphi_m = f^m s$ so that $\varphi_2$ is constant
by Lemma~\ref{lem2022-1-26-1}.
When $(M, g)$ has positive scalar curvature, we can see that, for $m=0, 1$, the function $\varphi_m$ attains its maximum and minimum at the points $x_0$ and $x_1$, respectively, 
where $f(x_0) = \min_M f$,  and $f(x_1) = \max_M f$. For $m \ge 3$, the function $\varphi_m$ 
attains its maximum and minimum at the points $x_1$ and $x_0$, respectively.
In fact, note that
$$
d\varphi_m = m\varphi_{m-1}df + f^m ds
$$
and, from  $fds = - 2s df$, we have
$$
d\varphi_m = (m-2) f^{m-1} s df
$$
and
\bea
Dd\varphi_m = (m-2)(m-3) f^{m-2}s df \otimes  df + (m-2) f^{m-1} s Ddf.\label{eqn2021-11-11-30}
\eea

We are ready to prove one of our main result.

\begin{thm}\label{prop2021-11-11-33}
Let $(M, g, f)$ be a compact Einstein-type manifold  with positive scalar curvature satisfying
(\ref{eq1}) for  constant $h$. 
If $\min_M s \ge nh$, then $(M,g)$ is Einstein.
\end{thm}
\begin{proof}
First of all, note that we must have $h \ge 0$. In fact, if $h <0$, then  $f<0$ on $M$  by Proposition~\ref{eqn2021-11-1-1}. So, letting $\min_M f = f(x_0)$, we have
$$
0 \le \Delta f(x_0) = f(x_0) s(x_0) - nh <0,
$$
which is a contradiction. Since it is reduced to Theorem 1.1 when $h=0$, we may assume that $h>0$ and so   the potential function $f$ is positive on $M$ by Proposition~\ref{eqn2021-11-1-1} again.
Since 
$$
\Delta f = fs - n h \ge (f-1)nh,
$$ 
considering the maximum point $x_1$ of $f$, we obtain
$$
f(x_1) \le 1,\quad \mbox{i.e.,}\quad f \le 1 \,\,\,\, \mbox{on $M$}.
$$
First, assume that 
$$
\max_M f = f(x_1) = 1.
$$
Since the function $fs$ attains its minimum at the point $x_1$, 
where $f(x_1) = \max_M f$,  we have
$$
0 \ge \Delta f(x_1) = (fs)(x_1) - n h.
$$
Since this implies $s(x_1) \le n h \le s(x_1)$, we have
$$
s(x_1) = n h.
$$
Now since $f^2 s = k$, constant, we have
$$
k = f(x_1)^2 s(x_1) = n h,
$$
and so
$$
fs = \frac{nh}{f}\ge nh.
$$
Therefore, $f$ must be constant since it is a subharmonic function.

Now, assume that $\max_M f = f(x_1) <1$. Then it is easy to compute that
$$
\Delta \ln(1-f) = \frac{|\n f|^2}{(1-f)^2} - \frac{fs-nh}{1-f},
$$
which can be written in the  following form
\be
\Delta \ln(1-f) - |\n \ln(1-f)|^2  - a \ln(1-f) = - a \ln(1-f) - \frac{fs}{1-f} 
+ \frac{nh}{1-f}, \label{eqn2021-11-12-10}
\ee
where $a$ is chosen to be a positive constant so that
\be
-a(1-f(x_0)\ln (1-f(x_0)) - (fs)(x_0) + n h >0.\label{eqn2021-11-13-1}
\ee
Since $\n f(\ln(1-f)) = -\frac{|\n f|^2}{1-f)}$ and $a>0$,  
the function $-a\ln(1-f)$ is non-decreasing in the $\n f$-direction. Also, since
$$
\n f\left( \frac{nh}{1-f}\right) =  \frac{nh|\n f|^2}{(1-f)^2}\quad\mbox{and}
\quad
\n f\left(\frac{f}{1-f}\right) = \frac{|\n f|^2}{1-f} + \frac{f|\n f|^2}{(1-f)^2},
$$
both two functions $ \frac{nh}{1-f}$ and  $ \frac{f}{1-f}$ are non-decreasing in the $\n f$-direction.
Since $s$ is non-increasing in the $\n f$-direction, the function $-\frac{fs}{1-f}$ is non-decreasing in the $\n f$-direction. Thus in view of (\ref{eqn2021-11-13-1}), we can see that the right-hand side
of (\ref{eqn2021-11-12-10}) is positive. Applying the maximum principle to
 (\ref{eqn2021-11-12-10}), we can conclude that $f$ must be constant.
\end{proof}

\begin{cor}
Let $(M, g, f)$ be a compact Einstein-type manifold  with positive scalar curvature satisfying
(\ref{eq1}) for a   non-zero constant $h$. 
If $\min_M s \ge nh$,  then, up to finite cover and scaling, $(M, g)$ is isometric 
to a standard sphere ${\Bbb S}^n$.
\end{cor}
\begin{proof}
Since $(M, g)$ is Einstein, it has positive Ricci curvature. By Myers' theorem, the fundamental group
$\pi_1(M)$ of $M$ is finite. Thus, up to finite cover and scaling, $(M, g)$ is isometric 
to a standard sphere ${\Bbb S}^n$.
\end{proof}

\section{Einstein-type manifolds with nonpositive Ricci curvature}

In this section, we will handle Einstein-type manifolds with nonpositive Ricci curvature.
First of all, we show that for such an Einstein-type manifold with nonpositive Ricci curvature, 
the square norm of the gradient of potential function satisfies the maximum principle.

\begin{lem} \label{lem930}
Let $(M^n,g,f, h)$ be an Einstein-type manifold with nonpositive Ricci curvature satisfying (\ref{eq1}) 
for a constant $h$. Then $|\nabla f|^2$ cannot attain its maximum in the interior. 
\end{lem}
\begin{proof}
By the assumption on the Ricci curvature, it is obvious that  $s-r(N,N)\leq 0,$ 
where $N=\nabla f/|\nabla f|$.
Since $\Delta f = fs$ and $f \n s = - 2 s \n f$, we have
$$
 \langle d\tr f, \nabla f\rangle = f \langle  \nabla s, \nabla f \rangle +s|\nabla f|^2=-s|\nabla f|^2.
$$
So, from the Bochner-Weitzenb\"ock formula, we obtain
$$
\frac 12 \tr |\nabla f|^2 +(s-r(N,N))|\nabla f|^2=|Ddf|^2\geq 0.
$$
Applying the maximum principle, we have the conclusion.
\end{proof}

Now, we present the proof of Theorem~\ref{thm930}.

\begin{thm} \label{thm930-1} 
Let $(M^n,g,f, h)$ be a complete noncompact Einstein-type manifold with nonpositive Ricci curvature satisfying (\ref{eq1}) with constant $h$.
If $f$ satisfies 
$$
\int_{M}|\nabla f|^2 < \infty,
$$
then $f$ is constant and $(M, g)$ is Einstein.
\end{thm}
\begin{proof}
For a cut-off function $\varphi$ (which will be determined later), we have
\bea
\frac 12 \int_M \varphi^2 \tr |\nabla f|^2&=&-2\int_M\varphi |\nabla f | \langle \nabla \varphi, \nabla |\nabla f|\rangle\\ 
&\leq & \int_M |\nabla f|^2 |\nabla \varphi|^2 +\varphi^2  |\nabla |\nabla f||^2.
\eea
Here we omit the volume form $dv_g$ determined by the metric $g$. From the Bochner formula derived in the proof of Lemma~\ref{lem930} with $N = \n f/|\n f|$, we have
$$
\frac 12\int_M \varphi^2 \tr |\nabla f|^2= \int_M \varphi^2 |Ddf|^2-(s-r(N,N))|\nabla f|^2 \varphi^2 .
$$
By combining these,
we obtain
\be
\int_M \left(|Ddf|^2 - |\n |\n f||^2\right)\vp^2  - \int_M (s-r(N,N))|\nabla f|^2 \varphi^2 \leq \int_M |\nabla f|^2 |\nabla \varphi|^2.
 \label{eqn2022-1-26-3}
 \ee
Now, take a geodesic ball $B_p(r)$ for some fixed point $p\in M$ and choose a cut-off function
$\varphi$ so that 
$$
{\rm supp} \, \varphi \subset B_p(r), \quad \varphi\vert_{B_p(\frac r2)}\equiv 1 , \quad 0\leq \varphi \leq 1,\quad  |\nabla \varphi|\leq \frac 1r.
$$ 
Substituting this into (\ref{eqn2022-1-26-3}) and using the  Kato's inequality, $ |\nabla |\nabla f||^2\leq  |Ddf|^2$, we have
$$   0\leq  \int_M \left(|Ddf|^2 - |\n |\n f||^2\right)  -\int_{B_p(\frac r2)} (s-r(N,N))|\nabla f|^2  \leq \frac 4{r^2}\int_{B_p(r)} |\nabla f|^2 .
$$
Letting $r\to \infty$, we obtain
$$ 
|Ddf|^2 = |\n |\n f||^2\quad \mbox{and}\quad (s-r(N,N))|\nabla f|^2=0
$$
on $M$. Hence, in particular, we have
$$ \frac 12 \tr |\nabla f|^2 =|\nabla |\nabla f||^2.$$ 
From
$$\frac 12 \tr |\nabla f|^2= |\nabla f|\tr |\nabla f|+|\nabla |\nabla f||^2,$$
we may conclude that $$|\nabla f|\tr |\nabla f|=0.$$
Again, for a cut-off function $\varphi$,
$$
0= \int_M \varphi^2 |\nabla f|\tr |\nabla f|= -2\int_M \varphi |\nabla f|\langle \nabla \varphi, \nabla |\nabla f|\rangle -\int_M \varphi ^2 |\nabla |\nabla f||^2.$$
By Young's inequality, we have
$$\int_M \varphi^2 |\nabla |\nabla f||^2 \leq \frac 12 \int_M \varphi^2 |\nabla |\nabla f||^2 + 2\int_M |\nabla \varphi|^2 |\nabla f|^2,$$
implying that
$$ \int_M \varphi^2 |\nabla |\nabla f||^2 \leq 4 \int_M |\nabla f|^2|\nabla \varphi|^2.$$
As shown above with the same cut-off function, the right-hand side tends to $0$. As a result, we have
$$ |\nabla |\nabla f||^2=0$$
on $M$, implying  that $Ddf\equiv 0$. Therefore, $g$ is Einstein with constant $f$.
\end{proof}

Before closing this section, we give an example of Einstein-type manifolds including space forms of 
nonnegative sectional curvature.

\begin{exam}[warped product]
{\rm
Let $g = dt^2 + \vp^2(t) g_0$ be a warped product metric on 
$M:= {\Bbb S}^1 \times_\vp {\Bbb S}^{n-1}$
or $M:= {\Bbb R} \times_\vp {\Bbb S}^{n-1}$, where $g_0$ is the standard round metric on ${\Bbb S}^{n-1}$.
Let $f = f(t)$ is a function defined on ${\Bbb S}^1$ or ${\Bbb R}$. It is easy to compute (cf. \cite{Be}) that
$$
{\rm Ric}_g = -(n-1)\frac{\vp''}{\vp}dt^2 - \left[\vp \vp''+(n-2)\vp'^2\right]g_0 + {\rm Ric}_{g_0}.
$$
It is also easy to see that the Hessian of $f$ with respect to the metric $g$ is given by
$$
Ddf(\partial_t, \partial_t)=f'',\quad
Ddf_{T{\Bbb S}^{n-1}} =\left(\frac {\vp'}{\vp}\right)f'g_{_{T{\Bbb S}^{n-1}}}, \quad 
Ddf(\partial_t, X)=0
$$
 for $X$ tangent to $T{\Bbb S}^{n-1}$. Here, $\partial_t = \frac{\partial}{\partial t}$ and
  `` $^{\prime}$ '' denotes the derivative taken with respect to $t\in {\Bbb S}$ or $t \in {\Bbb R}$.
 Since
 $$
 g_{_{T{\Bbb S}^{n-1}}} = \vp^2 g_0\quad \mbox{and}\quad {\rm Ric}_{g_0} = (n-2)g_0,
 $$
 it is easy to see that (\ref{eq1}) with $f$ and constant $h$  is equivalent to the following: 
\be
\left\{\begin{array}{ll}
f'' + (n-1)\frac{f\vp''}{\vp} + h = 0\\
f\left[-\vp \vp'' + (n-2)(1-\vp'^2)\right] - \vp \vp' f' - h \vp^2 = 0.
\end{array}\right.\label{eqn2022-1-26-6-1}
\ee
\begin{itemize}
\item[(I)] $\vp = c$ (constant)

This case does not happen, which means there are no solutions. 
In fact, if $\vp = c$ is constant, the first and second equations in (\ref{eqn2022-1-26-6-1}) are reduced to  $f''+h =0$ and 
$$
(n-2)f - hc^2 = 0,
$$
respectively. So, $f$ must be vanishing, which is  a contradiction.
\item[(II)] $\vp$ is non-constant

When $f>0$ on $M$, the first equation  in (\ref{eqn2022-1-26-6-1}) can be written as
\be
\frac{f''}{f} + \frac{h}{f} = - (n-1)\frac{\vp''}{\vp} =:\lambda.\label{eqn2022-1-26-5}
\ee
We assume that $\lambda$ is constant.
\begin{itemize}
\item[(i)] $\lambda >0$

Since the warping function $\vp$ must satisfy $\vp(0) = 0$ and $\vp'(0) = 1$,  the function $\vp$ satisfying (\ref{eqn2022-1-26-5}) has of the form
\be
\vp(t) = \sqrt{\frac{n-1}{\lambda}} \sin \sqrt{\frac{\lambda}{n-1}} t
\label{eqn2022-1-26-7}
\ee
and
\bea
f'' - \lambda f + h = 0.\label{eqn2022-1-26-8}
\eea
Note that $\vp$ is only working on a compact manifold $M= {\Bbb S}^1 \times_\vp {\Bbb S}^{n-1}$.
Next, note that the particular solution to this is given by
\be
f(t) = \frac{h}{\lambda} \,\,\, \mbox{(constant)} \label{eqn2022-1-26-9}
\ee
and a general solution is
$$
f(t) = ae^{\sqrt{\lambda}t} + \frac{h}{\lambda}.
$$
In particular, in case of compact manifold $M= {\Bbb S}^1 \times_\vp {\Bbb S}^{n-1}$, the function $f$ must be periodic, and so  $a=0$.
In conclusion, the compact manifold  $M= {\Bbb S}^1 \times_\vp {\Bbb S}^{n-1}$ with
$\vp$ and $f$ as in (\ref{eqn2022-1-26-7}) and (\ref{eqn2022-1-26-9}) is just an $n$-dimensional sphere ${\Bbb S}^n$ after a suitable scaling.
\item[(ii)] $\lambda =0$

In this case,  it is easy to see that 
$$
\vp(t) = t\quad \mbox{and}\quad f(t) = \mbox{constant}
$$
with $h = 0$ on ${\Bbb R}\times_\vp {\Bbb S}^{n-1}$  is the only possible case, and this case is just the flat Euclidean space ${\Bbb R}^n$.
\item[(iii)] $\lambda <0$

Letting $\mu = - \lambda >0$, as in (i) above, we have
$$
\vp(t) = \frac{1}{2}\sqrt{\frac{n-1}{\mu}} \left(e^{\sqrt{\frac{\mu}{n-1}}t} - e^{-\sqrt{\frac{\mu}{n-1}}t}\right) = \sqrt{\frac{n-1}{\mu}} \sinh \sqrt{\frac{\mu}{n-1}} t
$$
and
$$
f(t) = - \frac{h}{\mu} = \frac{h}{\lambda}.
$$
Thus, we can see this case is just the hyperbolic manifold.
\end{itemize}
\end{itemize}

}
\end{exam}

\section{Conformally flat Einstein-type manifolds}

Recall that, for an Einstein-type manifold satisfying (\ref{eq1}), if both $s$ and $h$ are constants and $f$ is not constant, then the scalar curvature $s$ must be zero.
In this section, we consider Einstein-type manifolds $(M, g, f, h)$ with zero scalar curvature satisfying (\ref{eq1}) with constant $h$ which $(M, g)$ is locally conformally flat.
Let us begin with the definition of Cotton tensor of a Riemannian manifold.
The Cotton tensor $C$ of a Riemannian manifold $(M, g)$ is defined by 
$$
C = d^D\left(r- \frac{s}{2(n-1)} g\right).
$$
Note that for a symmetric $2$-tensor $\xi$, $d^D\xi$ is defined as $d^D \xi(X, Y, Z) = D_X\xi(Y, Z) - D_Y\xi(X, Z)$for any vector fields $X, Y, Z$. It is well-known (cf. \cite{Be}) that the Cotton tensor
has a relation with the Weyl tensor ${\mathcal W}$ as follows
$$
 \mbox{div}\, {\mathcal W}=\frac {n-3}{n-2}C.
$$
It is also well-known (cf. \cite{Be}) that  if $\dim (M) = 3$, then $(M, g)$ is (locally) conformally flat if and only if the Cotton tensor $C$ is vanishing, and for $n \ge 4$, $(M, g)$ is (locally) conformally flat if and only if the Weyl tensor $\mathcal W$ is vanishing

When an Einstein-type manifold $(M, g, f, h)$ satisfying (\ref{eq1}) with constant $h$ is (locally)
conformally flat or  has zero Cotton tensor which is a little weak condition,  we can show that
$r(\n f, X) = 0$ for any vector field $X$ which is orthogonal to $\n f$, and this implies that every geometric data on $M$ is constant along each level hypersurface given by $f$, and the metric $g$ can be written as a warped product on the critical free set of $f$.
To do this, we introduce a $3$-tensor $T$ for Einstein-type manifolds defined as
$$ 
T= \frac 1{(n-1)(n-2)}\, i_{\nabla f} r \wedge g -\frac{s}{(n-1)(n-2)}df \wedge g + 
\frac{1}{n-2} df\wedge r,
$$
where $i_{\n f}$ denotes the interior product to the first component so that $i_{\n f}r(X) = r(\n f, X)$, and
$df \wedge \xi$ is defined as
$df\wedge \xi(X, Y, Z) = df(X)\xi(Y, Z) - df(Y)\xi(X, Z)$
for a symmetric $2$-tensor $\xi$ and vector fields $X, Y, Z$.
 
\begin{lem}\label{lem2021-1-11-2}
Let $(g,f,h)$ be a solution of $(\ref{eq1})$. Then
$$ f\, C=\tilde{i}_{\nabla f}{\mathcal W}-(n-1)T. $$
Here, $\tilde{i}_X$ is the interior product to the final factor by
$\tilde{i}_{\n f}  {\mathcal W} (X,Y,Z)= {\mathcal W}(X,Y, Z, \n f)$
for vector fields $X, Y, Z$.  
\end{lem}
\begin{proof}
cf. \cite{hy}.
\end{proof}

For an Einstein-type manifold $(M, g, f, h)$, let us denote $N=\nabla f/|\nabla f|$ and 
$\alpha: =r(N,N)$. It is clear that $\a$ is only well-defined on the set $\nabla f\neq 0$. 
However, since $|\a| \le |r|$, the function $\a$ can be defined on the whole $M$ as a $C^0$-function.

\begin{lem}\label{lem129} 
Let $(M,g,f,h)$ be a locally conformally flat Einstein-type manifold with constant $h$. Then we have $r(X,\nabla f)=0$ for $X$ orthogonal to $\nabla f$. In particular, the following holds.
\begin{enumerate}
\item[(1)] $|\nabla f|$ is constant on each level set of $f$.
\item[(2)]  $\alpha$ is constant on each level set of $f$. 
\item[(3)] Furthermore, if $s = 0$, then each level set of $f$ has constant sectional curvature.
\end{enumerate}
\end{lem}
\begin{proof}  By Lemma~\ref{lem2021-1-11-2} together with our assumptions, we have $T = 0$.
Let $\{e_i\}_{i=1}^n$  be an orthonormal frame  with $e_1=N$. 
Then, for $i\geq 2$, we have
$$ 0=(n-2)T(\nabla f, e_i, \nabla f)= \frac {n-2}{n-1}|\nabla f|^2r(e_i,\nabla f),
$$
which shows that  $r(X,\nabla f)=0$ for $X$ orthogonal to $\nabla f$.
\begin{itemize}
\item[(1)]
The property  $ r(e_i, \nabla f)=0$   for $i\geq 2$ also implies 
that $|\nabla f|$ is constant on each level set of $f$ by (\ref{eq1}). 

\item[(2)]
First of all,  note that  $D_NN=0$. In fact, we have 
$\langle D_NN, N\rangle =0$ trivially. 
Since
\bea
D_NN&=&N\left( \frac 1{|\nabla f|}\right)\n f+\frac 1{|\nabla f|}D_Ndf= -\frac 1{|\nabla f|}Ddf(N,N)N+\frac 1{|\nabla f|}D_Ndf\\
&=& -\frac 1{|\nabla f|} (f\alpha -h)N+\frac 1{|\nabla f|} \left[f r(N,\cdot)-h g(N, \cdot)\right],
\eea
we have  $ \langle D_NN, e_i\rangle =0$for $i\geq 2$.
Next, since $\langle D_{e_i}N, N\rangle =0$ and $C=0$, we have, for $i\geq 2$,
\bea
e_i(\alpha) &=&
D_{e_i}r(N,N) + 2 r(D_{e_i}N, N)\\
&=&
 D_N r(e_i, N)= - r(D_Ne_i, N)- r(e_i, D_NN)=-\langle D_Ne_i, N)\alpha =0.
\eea
\item[(3)]
From the fact $r(X,\nabla f)=0$ for $X$ orthogonal to $\nabla f$, we can write
 \be 
 i_{\nabla f}r =\alpha df \label{eq8}
 \ee
 as a $1$-form. Since   $T=0$ and $s=0$, we have
$$
df\wedge \left(r+ \frac{\a}{n-1}g\right) = 0,
$$
which shows
\be 
r_{ij}= -\frac {\alpha}{n-1}\d_{ij}\quad \mbox{and}\quad |r|^2=\frac n{n-1}\alpha^2.
\label{1230}
\ee
Therefore, the second fundamental form ${\rm II}$ on each level set $f^{-1}(c)$ of $f$ for a regular value $c$
 is given by 
\be 
{\rm II}(e_i, e_i)=\frac 1{|\nabla f|}Ddf(e_i, e_i)
= -\frac 1{|\nabla f|}\left(\frac {f \alpha}{n-1} +h\right)\label{eq11}
\ee
for $i\geq 2$.  Also, for $i,j\geq 2$ with $i\neq j$, we have
$$ R_{ijij}= \frac 1{n-2} (r_{ii}+r_{jj})= -\frac {2\alpha}{(n-1)(n-2)} .$$
Hence, by the Gauss equation with (\ref{eq11}), we conclude that $f^{-1}(c)$ has constant sectional curvature.
\end{itemize}
\end{proof}

We can see that, by equations (\ref{eq8}) and (\ref{1230}),  $\a$ is an eigenvalue of the Ricci curvature tensor $r$ with multiplicity one whose corresponding eigenvector is $\n f$, and $-\frac{\a}{n-1}$ is another eigenvalue with multiplicity $n-1$.

Even without the condition $C=0$ nor $\mathcal W = 0$, we can show that the scalar curvature $s$ is constant along each level hypersurface when $h$ is constant. In fact, from (\ref{eq83}), we have
$$
fds(e_i) = 0
$$
for $i \ge 2$.
Moreover, when $s = 0$,  note that,  by taking $d^D$ to (\ref{eq1}) 
$$
\tilde{i}_{\nabla f}R= d^DDdf=df\wedge r 
$$
since $d^Dr=0$.
In particular,  for $i,j\geq 2$ we obtain
\bea
 R_{1i1j}=-r_{ij} \label{eq353}
 \eea
with 
$$ R_{1i1i}= \frac {\alpha}{n-1}.
$$

\vspace{.2in}
From now on, we assume that $(M, g, f, h)$ is an Einstein-type manifold with constant $h$ 
and zero scalar curvature, $s=0$, and $f$ is proper so that each level hypersurface of $f$ is compact.
Let us denote by ${\rm Crit}(f)$ the set of all critical points of $f$.

\begin{lem}\label{lem2022-2-17-1}
Let $(M^n, g, f, h)$ be an Einstein-type manifold satisfying (\ref{eq1}) with constant $h$ and zero scalar curvature. Assume that $g$ is locally conformally flat. 
Then, on the set $M \setminus {\rm Crit}(f)$, we have
\be
N(\a) = \frac{n\a}{|\n f|}\left[\frac{f\a}{n-1}+h\right]\label{eqn2022-2-12-2}
\ee
and
\be
 NN(\a) = \frac{n}{n-1}\a^2 + \frac{n(n+1)\a}{|\n f|^2} \left[\frac{f\a}{n-1} + h\right]^2.\label{eqn2022-2-20-5}
 \ee
\end{lem}
\begin{proof}
Since $i_{\nabla f}r= \alpha df$ and $s=0$, we have
$$
 \langle \nabla \alpha, \nabla f\rangle -nh\alpha =\mbox{div}(i_{\nabla f}r)= f|r|^2
 $$
 and so
  \bea
 N(\a) = \frac{n\a}{|\n f|}\left[\frac{f\a}{n-1}+h\right].
 \eea
Taking the $N$-derivative of this and using  $N(|\n f|) = Ddf(N, N) = f\a - h$,  we obtain
\bea
 (f\a-h)N(\a) + |\n f|NN(\a)  = nh N(\a) + \frac{n}{n-1}|\n f|\a^2 + \frac{2nf\a}{n-1}N(\a).
 \eea
 So,
 \bea
 NN(\a) = \frac{n}{n-1}\a^2 + \frac{n(n+1)\a}{|\n f|^2} \left[\frac{f\a}{n-1} + h\right]^2.\label{eqn2022-1-8-4}
 \eea
 \end{proof}

 Denote by $\Sigma_c$ a connected component of the   level hypersurface $f^{-1}(c)$ of $f$ 
for a regular value $c\in {\mathbb R}$ which is assumed to be compact.

\begin{lem} \label{lem0107} 
Let $\Omega$ be a connected component of $M \setminus {\rm Crit}(f)$.
Then $(\Omega, g)$  is isometric to a warped product $(I \times \Sigma, dt^2+b^2g_{_\Sigma})$, 
where $I$ is an interval, $\Sigma$ is a connected component of the level hypersurface of $f$ contained in $\Omega$, and $g_{_\Sigma}$ is the induced metric on $\Sigma$ by $g$.  Moreover $g_{_\Sigma}$ is Einstein and $b^{n-1}\frac{d^2}{dt^2}b$ is constant on $I \times \Sigma$.
\end{lem}
\begin{proof} 
The second fundamental form ${\rm II}$ on $\Sigma$ is given by
\bea
 {\rm II}_{ij}= -\frac 1{|\nabla f|} \left( \frac {\alpha f}{n-1}+h\right)\delta_{ij}, \label{eq206}
\eea
which depends only on $f$ due to Lemma~\ref{lem129}. Thus, the metric can be written as 
$$ 
g=\frac {df}{|d f|}\otimes \frac {df}{|d f|} +b^2g_{_\Sigma},
$$
where $b=b(f)$ is a function of $f$ and $g_{_\Sigma}$ is the induced metric on $\Sigma$. 
Taking $dt=df/|\nabla f|$, we have $g=dt^2+b^2g_{_\Sigma}$.  In particular, it is easy to see that
$$
Ddf(\partial_t, \partial_t)=f'',\quad
Ddf_{T\Sigma}=\left(\frac {b'}{b}\right)f'g_{_{T\Sigma}}, \quad 
Ddf(\partial_t, X)=0
$$
 for $X$ tangent to $T\Sigma$. Here, ``\, $^{\prime}$\, '' denotes the derivative taken with respect to $t\in I$ so that $f'  = |\n f|$. For more details, we can refer \cite{lea},  \cite{cc} or \cite{MT}.
 
 By the standard calculation of the warped product and (\ref{eq1}) with $s=0$, we see that
 ${\rm Ric}(g_{_\Sigma})$ is given by 
\bea 
{\rm Ric}(g_{_\Sigma})
&=&
{\rm Ric}(g)\vert_{_{T\Sigma}}+\left((n-2)\left(\frac {b'}b\right)^2+\frac {b''}{b}\right)g_{_{\Sigma}}  \\
&=&\left(\frac {b'f'}{bf}+\frac h{f}+(n-2)\left(\frac {b'}b\right)^2+\frac {b''}{b}\right)g_{_{\Sigma}},
\eea
which implies that $g_{_\Sigma}$ is Einstein.

Next, since $D_NN=0$, we have $Ddf = f''dt^2 + bb'f' g_{_\Sigma}$ and so
\be
\Delta f = f'' + (n-1)\frac{b'}{b}f'.\label{eqn2022-2-5-1}
\ee
It is easy to compute the Ricci curvature of $g$ which is given by
$$
{\rm Ric}_g = -(n-1)\frac{b''}{b}dt^2 - \left[b b''+(n-2)b'^2\right]g_{_\Sigma} + {\rm Ric}_{g_{_\Sigma}},
$$
So, 
\bea
{\rm Ric}_g(\partial_t, \partial_t) = -(n-1)\frac{b''}{b}\label{eqn2022-2-3-2}
\eea
and the Einstein-type equation (\ref{eq1}) with $f$ is equivalent to the following: 
\be
\left\{\begin{array}{ll}
f'' + (n-1)\frac{ b''}{b}f + h = 0\\
f\left[-b b'' + (n-2)(\kappa_0-b'^2)\right] - b b' f' - h b^2 = 0.
\end{array}\right.\label{eqn2022-1-26-6}
\ee
Here ${\rm Ric}_{g_{_\Sigma}} = (n-2)\kappa_0 g_{_\Sigma}.$ 
Recall that we assume $h$ is constant. 
On the other hand, since $\Delta f = -nh$, we have
\be 
 -(\Delta f)g +Ddf-fr= (n-1)h g. \label{vst1}
\ee
Thus, by substituting  the pair $(\partial_t, \partial_t)$ into (\ref{vst1}),  we obtain
\be 
fb''- f' b'=hb. \label{eq202}
\ee
Taking the derivative of (\ref{eq202}) and substituting the first identity in (\ref{eqn2022-1-26-6}), 
we obtain
\bea
 bb''' + (n-1)b'b'' = 0\label{eqn2022-1-28-2}
\eea
and so $(b^{n-1}b'')' = 0.$
\end{proof}

In Lemma~\ref{lem0107}, if $b$ is constant so that $g$ is a product metric, then,  obviously, $\Sigma$ is flat and $(M, g)$ is Ricci-flat. 
We may also assume that the interval $I$ is given by $[t_0, t_1)$ and $f$ has a critical point in $f^{-1}(t_1)$ when ${\rm Crif}(f) \ne \emptyset$.
Also, by Lemma~\ref{lem0107}, we can let
\be
b^{n-1}b'' = a_0\label{eqn2022-1-28-3}
\ee
for some constant $a_0$.
Moreover, by warped product formula from $g = dt^2 + b^2 g_{_\Sigma}$, we have 
\be
\a = {\rm Ric}(\partial_t, \partial_t)=- (n-1)\frac {b''}b= -(n-1)\frac{a_0}{b^n}.\label{eqn2022-2-21-10-1}
\ee
That is,
\be 
\frac {b''}b=-\frac {\alpha}{n-1}= \frac {a_0}{b^n}. \label{eq206-1}
\ee
So, since $\a b^n = -(n-1)a_0$, we have
\be
 \a' b + n \a b' =0.\label{eqn2022-2-12-1-1}
\ee

From now on, we will show $a_0$ must be zero, which implies $(M, g)$ is Ricci-flat.
First, if we assume $a_0 \ne 0$, we can show that the warping function $b$ can be extended beyond the critical point of $f$ and so the warped product metric is still valid beyond the critical points of $f$.

\begin{lem}\label{lem2022-2-26-10-100}
Under the hypotheses of Lemma~\ref{lem0107}, assume that $I =[t_0, t_1)$ so that $f$ has a critical point at $t=t_1$.
If $a_0 \ne 0$, then the warping function $b$ can be extended smoothly beyond $t_1$ satisfying (\ref{eqn2022-1-28-3}). Moreover, $\Sigma_1:=f^{-1}(t_1)$ is a smooth hypersurface, and
$f\a +(n-1) h = 0$ on the set $\Sigma_1$.
\end{lem}
\begin{proof}
Note that the scalar curvature of the metric $g$ is given by
$$
0=s=-2(n-1)\frac {b''}b +\frac {s_c}{b^2}-(n-1)(n-2)\left(\frac {b'}b\right)^2 
$$
with $s_\Sigma= (n-1)(n-2) \kappa$, where $s_\Sigma$ is the (normalized) scalar curvature of $\Sigma = f^{-1}(t_0)$. So, we obtain
\be
2\frac {b''}b +(n-2)\left(\frac {b'}b\right)^2 = (n-2)\frac {\kappa}{b^2}.
\label{eqn2022-2-3-1-1} 
\ee
This can be written in the following form
\be
(n-2)b'^2+2a_0b^{2-n} =(n-2)\kappa. \label{eqn2022-1-28-6-1}
\ee
If $a_0 \ne 0$, then (\ref{eqn2022-1-28-6-1}) shows
$$
\liminf_{t\to t_1} b(t) >0.
$$
In fact, if $a_0>0$ and $\dis{\liminf_{t\to t_1} b(t) =0}$, then we have
 (LHS) $\to \infty$ and (RHS) $=(n-2)\kappa$ in (\ref{eqn2022-1-28-6-1}), 
 a contradiction. 
Now assume $a_0<0$ and $\dis{\liminf_{t\to t_1} b(t) =0}$. Then, by (\ref{eqn2022-1-28-3}),
$$
\lim_{t\to t_1} \frac{b''}{b}(t) =  -\infty.
$$
Since $f$ is well-defined around $t = t_1$, by (\ref{eq202}), 
 we must have
$$
 \lim_{t\to t_1} \frac{f'b'}{b}(t) = -\infty.
 $$
However, this contradicts (\ref{eqn2022-2-5-1}) since $\Delta f$ is well-defined on $t=t_1$.

Suppose, now,  there exists $t_k\to t_1$ such that $b(t_k) \to \infty$ as $k \to \infty$. 
It follows from (\ref{eqn2022-1-28-6-1})  that
$$
\left(\frac{b'(t_k)}{b(t_k)}\right)^2 \to 0.
$$
So, by (\ref{eqn2022-1-28-3}) or (\ref{eqn2022-2-3-1-1}), we have
\bea
\lim_{k\to \infty} \frac{b''}{b} = \lim_{k\to \infty} \frac{a_0}{b^{n}} =0. \label{eqn2022-1-30-3}
\eea
Since $f'$ is well-defined on $M$, we have, by (\ref{eq202}),
$$
h = \lim_{k\to \infty} f\frac{b''}{b} -\lim_{k\to \infty} f'\frac{b'}{b} = 0,
$$
which  contradicts $h >0$. Therefore we have
$$
\limsup_{t\to t_1} b(t) < +\infty.
$$
It follows that $C^{-1} \le  b \le C$ on $I=[t_0, t_1)$ for some $C>0$.
In particular, $b$ can be extended smoothly beyond $t_1$ satisfying (\ref{eqn2022-1-28-3}), 
and $f^{-1}(t_1)$ is a smooth hypersurface.


Since $f$ has no critical points on $t_0\le t <t_1$, every level set $f^{-1}(t)$ for $t_ 0\le t <t_1$ is homotopically equivalent. Thus, at $t=t_1$, $f^{-1}(t_1)$ would be a hypersurface which is homotopically equivalent to $\Sigma$. Since $f^{-1}(t_1)$ is a hypersurface, the mean curvature $m$ is well-defined. Since, around $t = t_1$, we have
$$
-nh = \Delta f = Ddf(N, N) + m |\n H| = f\a - h + m |\n f|,
$$
i.e.,
\bea
f\a  +(n-1)h = -m |\n f|,\label{eqn2022-2-25-1-1}
\eea
by letting $t \to t_1$, we obtain
\bea
f\a + (n-1)h = 0\label{eqn2022-2-24-1-1}
\eea
on the set $f^{-1}(t_1)$. 
\end{proof}

\begin{lem}\label{lem2022-2-26-1}
Under the hypotheses of Lemma~\ref{lem0107} with $a_0 \ne 0$, 
let $\Sigma_1:=f^{-1}(t_1)$ be a smooth hypersurface. Assume  $f\a + (n-1)h = 0$ on the set $\Sigma_1$. 
Then, we have  the following.
\begin{itemize}
\item[$(1)$] $N(\a) = \a' = 0, \,\,\, NN(\a)=  \a'' >0$ on the set $\Sigma_1$.
\item[$(2)$] $\Sigma_1$ is totally geodesic.
\end{itemize}
\end{lem}
\begin{proof}
\begin{itemize}
\item[(1)]
 If $\Sigma_1=f^{-1}(t_1)$ contains a  critical point of $f$, then
 $$
 b' = \frac{db}{df}\cdot \frac{df}{dt} = 0
 $$ 
 at $t=t_1$, and so $\a' = 0$ at $t=t_1$ by (\ref{eqn2022-2-12-1-1}).
If $f$ has no critical point on $\Sigma_1$, by Lemma~\ref{lem2022-2-17-1}, we have $\a'=0$ at $t=t_1$.
The second inequality $\a'' >0$ follows from  (\ref{eqn2022-2-20-5}).

\item[(2)]
Let $\nu$ be the outward unit normal vector field on $\Sigma_1$ such that
$$
\lim_{t\to t_1-}N = \nu.
$$
Let $\{e_i\}$ be a local frame around $\Sigma_1$ such that $e_1 = N = \frac{\n f}{|\n f|}$ and then
$$
\lim_{t\to t_1-}e_1 = \nu.
$$
Then, on the set $f^{-1}(t_1-\e)\cap \Omega$ for a sufficiently small $\e>0$,  we have
$$
D_{e_i}N = - \frac{1}{|\n f|}\left(\frac{f\a}{n-1} +h\right) e_i = \frac{b'}{b}e_i.
$$
By letting $\e \to 0$, by (1) above, we obtain
$$
D_{e_i} \nu = 0,
$$
which shows $\Sigma_1$ is totally geodesic.

\end{itemize}
\end{proof}

\vspace{.12in}
\noindent
{\bf Proof of Theorem~\ref{thm129}.}\,\, 
Since both $h$ and $s$ are constants, we have
$sdf = 0$ by (\ref{eq83}) so that either $f$ is constant on $M$ or $s=0$. If $f$ is constant, $(M, g)$ is Einstein satisfying $fr = hg$. 
If  $M$ is compact (without boundary), then from $\Delta f = -nh$, $f$ must be constant 
and so $(M, g)$ is Einstein. 

Now, assume $(M, g, f, h)$ is a complete non-compact Einstein-type manifold with constant $h>0$ 
and zero scalar curvature, $s=0$. Moreover, we assume that $f$ is not constant and a proper map.
If  $f$ contains an isolated critical point, 
then $a_0$ must be zero  by Lemma~\ref{lem2022-2-26-10-100}. 
So, by  (\ref{1230}) and (\ref{eq206-1}), $(M, g)$ is Ricci-flat.
Assume that every  level set of $f$ is a compact smooth hypersurface.
There are two  cases; either $f$ has no critical points, or 
  $M$ contains a level hypersurface  given by $f$ consisting of critical points of $f$.
For the latter case, if $f^{-1}(t_1)$ is a smooth hypersurface consisting of critical points, and
we have $b'(t_1)= 0=\a'(t_1)$ and $\a''(t_1)>0$ as in the proof of Lemma~\ref{lem2022-2-26-1}.
In particular, $f^{-1}(t_1)$ is totally geodesic. Thus, by Lemma~\ref{lem0107}, the metric $g$ can be written as a warped product metric
$$
g = dt^2 + b(t)^2 g_{_\Sigma}
$$ 
globally on  $M=(-\infty, \infty) \times \Sigma$, $[a, \infty)\times \Sigma$, or $(-\infty, a]\times \Sigma$ 
for a hypersurface $\Sigma$ of $M$ with constant sectional curvature.   By (negative) parametrization, we only consider first two cases, and the second case corresponds to incomplete Einstein-type manifold.

\begin{itemize}

\item[(1)] $a_0 <0$

In this case, we have $\a >0$  and $b'' <0$ on $M$ by (\ref{eq206-1}).
 If $f^{-1}(t_0)$ is a smooth hypersurface consisting of critical points, then, by Lemma~\ref{lem2022-2-26-1}, it is totally geodesic, and $ f\a = -(n-1)h$, 
  $\a' = 0$ and $\a''>0$ at $t=t_0$. In particular, the function $\a$ attains its local minimum at $t = t_0$.
    Since $\a  b^n = -(n-1)a_0$, $b$ attains its local maximum at $t=t_0$, i.e., $b'(t_0) = 0$ and
 $b''(t_0) <0$. Moreover since
 $$
 \left(\frac{b'}{b}\right)' = \frac{b''}{b} - \frac{b'^2}{b^2} = - \frac{\a}{n-1} - \frac{b'^2}{b^2} <0,
 $$
 we have $ b' <0$ for $t >t_0$, and $b'>0$ on $t<t_0$. This shows that $b$ is convex upward and attains its global maximum at $t=t_0$.  However this is impossible since $M$ is complete and $b>0$.
 
 Now assume that 
 $$
 f\a >- (n-1)h
 $$
 on $M$ so that $\a' >0, \a'' >0$ and $f$ has no critical points. Then
 $$
 b'<0\quad \mbox{and}\quad b'' <0.
 $$
 This is also impossible since $b>0$ and $M$ is complete and noncompact.

Since we assumed $h>0$, the case that $f\a < -(n-1)h$ on $M$ does not happen
 by considering the set $f>0$. Remark that the last two cases correspond to the case
 $M= (-\infty, \infty)\times \Sigma$ by Lemma~\ref{lem2022-2-26-10-100}.


\item[(2)]  $a_0 >0$

In this case, we have $\a <0$  and $b'' >0$ on $M$ by (\ref{eqn2022-2-21-10-1}).
 If there exists a $t_1 \in {\Bbb R}$ such that $f\a = -(n-1)h$ at $t=t_1$, we have
  $b' =0$ at $t=t_1$ as above. Since $b>0$, we have
  $$
 \lim_{t\to \pm \infty} b(t) = \infty
 $$
  unless $b$ is constant. By (\ref{eq206-1}),
  \bea
  \lim_{t\to \pm \infty} \a(t) = 0.\label{eqn2022-2-28-1}
 \eea
 However, since $\a'(t_1)= 0, \a''(t_1)>0$ and $\a <0$ on $\Bbb R$,  this is impossible.
 
 Now, assume that 
 \be
 f\a >- (n-1)h\label{eqn2022-2-21-10}
 \ee
 on $M$ so that $\a' <0$.  Note that, from Lemma~\ref{lem2022-2-26-10-100}, this corresponds to the first case $M = (-\infty, \infty)\times \Sigma$ with no critical points of $f$. Then
 $$
 b'< 0 \quad \mbox{and}\quad b'' >0.
 $$
 This shows that
 $$
 \lim_{t\to - \infty} b(t) = \infty$$
 and so, by (\ref{eqn2022-2-21-10-1}),
 $$
 \lim_{t \to -\infty} \a = 0.
 $$
On the other hand, on the set $f>0$, we have 
$$
- \frac{(n-1)h}{f}< \a <0$$
by (\ref{eqn2022-2-21-10}), and so  we have
$$
 \lim_{t \to \infty} \a = 0 .
 $$
However, this is impossible since $\a<0$ and $\a' <0$ on $M$.
 
 Finally,  assume  that $f\a < -(n-1)h$ on $M$.
 Recall that, from Lemma 5.5, this also corresponds to the first case
$M = (-\infty, \infty) \times \Sigma$ with no critical points of $f$. As above, considering the set $f<0$, we can see that this is impossible because of $h>0$.
\end{itemize}
 Hence, we must have $a_0 = 0$ and so $(M, g)$ is Ricci-flat.
 \hfill $\Box$

\section{Final Remarks}

For Einstein-type manifolds satisfying (\ref{eq1}), even though the potential function $f$ has a relation to $h$ as $h = \frac{1}{n}(fs - \Delta f)$,  if $h$ is not constant, it is not easy to obtain rigidity results since the function $h$ causes  various situations. Nonetheless, if we give some
 constraints on the function $h$, we can have a little minor results. 
We say that $h$ has defined {\it weakly signal} if either $h\geq 0$ on $M$ or $h\leq0$ on $M$.
With this assumption, we can show that there are no critical points on $f^{-1}(0)$ unless $f^{-1}(0)$ is empty, and this implies the scalar curvature is vanishing on the set $f^{-1}(0)$ as in the proof of Lemma~\ref{lem83}.

\vskip .5pc
On the other hand, when $h=csf$ for some constant $c$, we have the following result under nonnegative Ricci curvature condition.

\begin{lem}
Let $(M^n,g,f, h)$ be a compact Einstein-type manifold with boundary $\partial M$ with $h=csf$.  If $(M,g)$ has nonnegative Ricci curvature with for $c< \frac 1{2(n-1)}$ or $c\geq \frac 1n$, then $|\nabla f|^2$ should have its maximum on $\partial M$.
\end{lem}
\begin{proof}
Note that $s\geq 0$. By (\ref{eq83})
and 
$$
\nabla h= c (s\nabla f+f\nabla s),
$$
we have
$$ 
(1-2nc+2c)f\nabla s=2(nc-c-1)s\nabla f,
$$
implying that
$$
d\tr f =(1-nc) (sdf+fds)= \frac {1-nc}{2nc-2c-1}sdf.
$$
By Bochner-Weitzenb\"ock formula and the assumption on the Ricci curvature, we have
$$ 
\frac 12 \tr |\nabla f|^2- \frac {nc-1}{2nc-2c-1}s|\nabla f|^2 =|Ddf|^2+r(\nabla f, \nabla f)\geq 0.
$$
From the condition on the constant $c$, we have
$$ 
\frac {nc-1}{2nc-2c-1} >0.
$$
Our Lemma follows from the maximum principle.
\end{proof}

\begin{cor} Let $h=csf$ for a constant $c$ satisfying  $c< \frac 1{2(n-1)}$ or $c\geq \frac 1n$. If $(M,g,f,h)$ is a closed Einstein-type manifold with nonnegative Ricci curvature,  then $(M,g)$ is Einstein.
\end{cor}

\vskip 0.3cm
\noindent
{\bf Acknowledgment:} The authors would like to express their gratitude to the referees for
valuable suggestions.
The first-named author was supported by the National Research
 Foundation of Korea (NRF-2019R1A2C1004948) and the second-named author was supported
  by the National Research Foundation of Korea (NRF-2018R1D1A1B05042186).


\bigskip
\noindent Gabjin Yun\\
Department of Mathematics and The Natural Science Research Institute\\
Myongji University \\
Myongji-ro 116, Cheoin-gu, Yongin, Gyeonggi-do, 17058, Korea \\
{\tt E-mail:gabjin@mju.ac.kr} \\

\bigskip
\noindent Seungsu Hwang (corresponding author)\\
{Department of Mathematics}\\
Chung-Ang University\\
84 HeukSeok-ro DongJak-gu, Seoul, 06974, Republic of Korea \\
{\tt E-mail:seungsu@cau.ac.kr}\\

\end{document}